\documentclass{article}
\usepackage[margin=1in]{geometry}
\usepackage{mathtools}
\usepackage{amssymb}
\usepackage{amsthm}

\newtheorem{theorem}{Theorem}
\newtheorem{lemma}[theorem]{Lemma}

\theoremstyle{definition}
\newtheorem{definition}{Definition}

\newcommand{\RR}{\mathbb{R}}

\DeclareMathOperator*{\st}{s.t.}
\newenvironment{myproof}{\begin{proof}}{\end{proof}}

\begin{document}

\title{On Upper Bounding Shannon Capacity of Graph Through Generalized Conic Programming}

\author{Yingjie Bi and Ao Tang \\ Cornell University, Ithaca, NY, 14850}
\date{}
\maketitle
\begin{abstract}
The Shannon capacity of a graph is an important graph invariant in information theory that is extremely difficult to compute. The Lov\'asz number, which is based on semidefinite programming relaxation, is a well-known upper bound for the Shannon capacity. To improve this upper bound, previous researches tried to generalize the Lov\'asz number using the ideas from the sum-of-squares optimization. In this paper, we consider the possibility of developing general conic programming upper bounds for the Shannon capacity, which include the previous attempts as special cases, and show that it is impossible to find better upper bounds for the Shannon capacity along this way.
\end{abstract}

\section{Introduction}

The \emph{Shannon capacity of a graph} is a graph invariant originated from computing the maximum achievable rate to transmit information with zero possibility of error through a noisy channel \cite{Shannon1956}. To state the definition of Shannon capacity, we need the following notions in graph theory: For an undirected graph $G$, let $V(G)$ and $E(G)$ be its vertex set and edge set, respectively. Let $\alpha(G)$ be the independence number (aka stability number) of $G$, i.e., the size of the maximum independent set in $G$. For two vertices $i,j \in V(G)$, the notation $i \sim_G j$ means either $i=j$ or $(i,j) \in E(G)$. The \emph{strong product} $G \boxtimes H$ of two graphs $G$ and $H$ is a graph such that
\begin{itemize}
\item its vertex set $V(G \boxtimes H)$ is the Cartesian product $V(G) \times V(H)$ and
\item $(i,j) \sim_{G \boxtimes H} (k,l)$ if and only if $i \sim_G k$ and $j \sim_H l$.
\end{itemize}
The Shannon capacity $\Theta(G)$ of graph $G$ is defined by
\[
\Theta(G)=\sup_k\sqrt[k]{\alpha(G^k)},
\]
where $G^k$ is the strong product of $G$ with itself for $k$ times.

The Shannon capacity is unknown for most graphs, including certain simple cases such as odd cycles $C_{2n+1}$ when $n \geq 3$. By definition, for any positive integer $k$, $\sqrt[k]{\alpha(G^k)}$ provides a direct lower bound for the Shannon capacity $\Theta(G)$, although it is still hard to calculate due to the NP-hardness of maximum independent set problem and the exponential growth of the size of $G^k$. Finding a good upper bound for $\Theta(G)$ is even more difficult. One well-known upper bound is the Lov\'asz number $\vartheta(G)$ proposed in \cite{Lovasz1979}, which can be efficiently computed by solving a semidefinite program (SDP). The most famous application of Lov\'asz number is the establishment of the Shannon capacity for the pentagon graph $C_5$:
\[
\sqrt 5=\sqrt{\alpha(C_5^2)} \leq \Theta(C_5) \leq \vartheta(C_5)=\sqrt 5.
\]
However, for 7-cycle $C_7$, $\vartheta(C_7) \approx 3.3177$, while the best known lower bound \cite{PS2019} at the time of writing is
\[
\Theta(C_7) \geq \sqrt[5]{\alpha(C_7^5)} \geq \sqrt[5]{367} \approx 3.2578.
\]
Determining the exact value for the Shannon capacity $\Theta(C_7)$ remains an open problem.

One interesting direction is to look for a tighter upper bound for the Shannon capacity than the Lov\'asz number. Since the definition of the Shannon capacity is closely related to the independence number, and in fact the Lov\'asz number itself can be derived from approximating the independence number of a graph, it is tempting to find better upper bounds for the Shannon capacity by using tighter approximations for the independence number. The major challenge here is to ensure that the new approximation is still an upper bound for the Shannon capacity. In Section~\ref{sec:conic}, we will look at general conic programming approximation for the independence number, which is a natural generalization of the SDP-based Lov\'asz number. Next, in Section~\ref{sec:product}, we will propose a condition called the \emph{product property} over the cones appeared in the above approximate optimization problem. This property guarantees that the optimal value of the approximation is an upper bound for the Shannon capacity. Surprisingly, in Section~\ref{sec:optimal} it is shown that the semidefinite cone used by the Lov\'asz number is the largest cone with such a property, thus ruling out the possibility of improving the estimation of the Shannon capacity along this way.

\section{Conic Programming for the Independence Number}\label{sec:conic}

In this section, we will first formulate the maximum independent set problem as a copositive program. If the semidefinite cone is used as an inner approximation for the copositive cone in this program, the obtained objective value is exactly the Lov\'asz number. As a generalization, we consider all the possible cones that are subsets of the copositive cone, and the corresponding conic programs will be the candidates to generate better upper bounds for the Shannon capacity. Before we start, we summarize the common notations used in the paper below:
\begin{itemize}
\item $\RR_+^n$ is the set of $n \times 1$ nonnegative column vectors.
\item $J_n$ is the $n \times n$ matrix of all ones.
\item $\mathcal S_n$ is the cone of $n \times n$ symmetric matrices.
\item $\mathcal P_n$ is the cone of $n \times n$ positive semidefinite matrices.
\item $\mathcal N_n$ is the cone of $n \times n$ nonnegative symmetric matrices.
\item $\mathcal C_n$ is the cone of $n \times n$ copositive matrices, i.e., all symmetric matrices $Q \in \mathcal S_n$ such that $x^TQx \geq 0$ for any $x \in \RR_+^n$.
\end{itemize}

Our starting point is the Motzkin-Straus theorem, which gives the exact value of the independence number of a graph:

\begin{theorem}[Motzkin-Straus]\label{thm:motzkin}
If $A$ is the adjacency matrix of a graph $G$ with $n$ vertices, then the independence number of $G$ is given by
\[
\frac{1}{\alpha(G)}=\min_{x \in \RR_+^n,\sum_ix_i=1}x^T(I+A)x.
\]
\end{theorem}

In \cite{KP2002}, the optimization problem in Theorem~\ref{thm:motzkin} is converted into the following equivalent form:
\begin{equation}\label{eq:copold}
\begin{split}
\alpha(G)=\min \quad & \lambda \\
\st \quad & \lambda(I+A)-J_n \in \mathcal C_n.
\end{split}
\end{equation}
In order to make the above problem \eqref{eq:copold} closer to the formulation for the Lov\'asz number, we are going to further rewrite it as follows:
\begin{equation}\label{eq:cop}
\begin{split}
\min \quad & \lambda \\
\st \quad & Y-J_n \in \mathcal C_n, \\
& Y_{ii}=\lambda, \quad \forall i=1,\dots,n, \\
& Y_{ij}=0, \quad \forall i \not\sim_G j, \\
& Y \in \mathcal S_n.
\end{split}
\end{equation}
Since problem \eqref{eq:copold} can be viewed as problem \eqref{eq:cop} with the additional constraint $Y=\lambda(I+A)$, problem \eqref{eq:cop} is a relaxation of the original problem \eqref{eq:copold}. To show that these two problems are indeed equivalent, the following property of copositive matrices will be useful:

\begin{lemma}\label{lem:copmat}
Assume $Q$ is a copositive matrix whose diagonal entries are all equal to $\mu$. $R$ is another symmetric matrix of the same size. If for each entry of $R$ either $R_{ij}=Q_{ij}$ or $R_{ij}=\mu$, then $R$ is also copositive.
\end{lemma}
\begin{myproof}
We only need to consider the case in which $R=Q$ except for some off-diagonal entry $R_{st}=\mu$ (and also $R_{ts}=\mu$), since the general result can be obtained by repeating the same argument for each difference between $R$ and $Q$. For any $x \in \RR_+^n$ with $\sum_ix_i=1$,
\begin{equation}\label{eq:bilinear}
x^TRx=\mu x_s^2+\mu x_t^2+2\mu x_sx_t+\sum_{\substack{(i,j) \neq (s,s),\\(s,t),(t,s),(t,t)}}Q_{ij}x_ix_j.
\end{equation}
Fix $x_i$, $i \neq s,t$, as constants and regard $x^TRx$ as a function of $x_s$ by replacing
\[
x_t=1-x_s-\sum_{i \neq s,t}x_i.
\]
Then the first part of \eqref{eq:bilinear}
\[
\mu x_s^2+\mu x_t^2+2\mu x_sx_t=\mu(x_s+x_t)^2=\mu\left(1-\sum_{i \neq s,t}x_i\right)^2
\]
becomes a constant. Since the remaining terms in \eqref{eq:bilinear} are all linear functions of $x_s$, $x^TRx$ is also linear as a function of $x_s$ and thus must achieve the minimum when $x_s=0$ or $x_s=1$. However, in both cases, $x^TRx=x^TQx \geq 0$, which implies that $R$ is also copositive.
\end{myproof}

Now we can prove that the problems \eqref{eq:copold} and \eqref{eq:cop} have the same optimal value. Consider an arbitrary feasible solution $(\lambda,Y)$ to problem \eqref{eq:cop}. Let
\[
Q=Y-J_n, \quad R=\lambda(I+A)-J_n.
\]
All the diagonal entries of $Q$ are $\lambda-1$. By Lemma~\ref{lem:copmat}, the matrix $R$ is also copositive and thus $\lambda \geq \alpha(G)$ by \eqref{eq:copold}. On the other hand, the solution $\lambda^*=\alpha(G)$, $Y^*=\alpha(G)(I+A)$ is feasible to \eqref{eq:cop}, so it must be optimal.

The copositive cone constraint in \eqref{eq:cop} makes the problem hard to solve. If we substitute the copositive cone $\mathcal C_n$ in \eqref{eq:cop} with the semidefinite cone $\mathcal P_n$, the optimal value for the modified problem is exactly the Lov\'asz number\footnote{The common definition of the Lov\'asz number that appears in the literature is the dual form of ours.}, which will be denoted as $\vartheta(G)$. Since $\mathcal P_n \subseteq \mathcal C_n$, we immediately get $\alpha(G) \leq \vartheta(G)$. Naturally, to find a tighter bound for the Shannon capacity $\Theta(G)$, we can replace the copositive cone $\mathcal C_n$ in \eqref{eq:cop} by some cone between $\mathcal C_n$ and $\mathcal P_n$, which may lead to some problem whose optimal value is potentially between the Shannon capacity $\Theta(G)$ and the Lov\'asz number $\vartheta(G)$.

The above discussion illuminates us to construct more general approximations for the independence number $\alpha(G)$ by introducing an arbitrary cone $\mathcal A_n \subseteq \mathcal C_n$ into the problem
\begin{equation}\label{eq:conic}
\begin{split}
\min \quad & \lambda \\
\st \quad & Y-J_n \in \mathcal A_n, \\
& Y_{ii}=\lambda, \quad \forall i=1,\dots,n, \\
& Y_{ij}=0, \quad \forall i \not\sim_G j, \\
& Y \in \mathcal S_n.
\end{split}
\end{equation}

In the case when the cone $\mathcal A_n$ is chosen to be the semidefinite cone $\mathcal P_n$, the above problem \eqref{eq:conic} gives the Lov\'asz number $\vartheta(G)$. To provide some other examples of $\mathcal A_n$, one can approximate the copositive cone $\mathcal C_n$ based on sum-of-squares programming \cite{KP2002,DR2010}. Note that the copositivity of a matrix $Q \in \mathcal S_n$ is equivalent to the condition
\begin{equation}\label{eq:sospoly}
p_Q(x)=\sum_{i,j}Q_{ij}x_i^2x_j^2 \geq 0, \quad \forall x \in \RR^n.
\end{equation}
Like determining copositivity, it is NP-hard to decide whether the polynomial $p_Q(x)$ is nonnegative or not. However, if $p_Q(x)$ can be written as a sum of squares, i.e.,
\[
p_Q(x)=\sum_kg_k^2(x),
\]
where $g_k(x)$ are arbitrary polynomials of $x \in \RR^n$, then clearly $p_Q(x)$ is nonnegative. All symmetric matrices $Q \in \mathcal S_n$ whose corresponding polynomial $p_Q(x)$ given by \eqref{eq:sospoly} is a sum of squares constitute a cone, which will be denoted as $\mathcal C_n^{(0)}$ in the following. By the above discussion, $\mathcal C_n^{(0)} \subseteq \mathcal C_n$, and furthermore it is tractable to determine whether a matrix $Q$ is in the cone $\mathcal C_n^{(0)}$ through SDP. In fact, $\mathcal C_n^{(0)}$ has a simple characterization \cite{Parrilo2000}:
\[
\mathcal C_n^{(0)}=\mathcal P_n+\mathcal N_n.
\]
In other words, the polynomial $p_Q(x)$ is a sum of squares if and only if the matrix $Q$ is a sum of a positive semidefinite matrix and a nonnegative symmetric matrix.

For any graph $G$, the optimal value of problem \eqref{eq:conic}, in which $\mathcal A_n=\mathcal C_n^{(0)}$, is called $\vartheta'(G)$, the Schrijver $\vartheta'$-function \cite{Schrijver1979}. Since
\[
\mathcal P_n \subseteq \mathcal C_n^{(0)} \subseteq \mathcal C_n,
\]
we have
\[
\alpha(G) \leq \vartheta'(G) \leq \vartheta(G).
\]
Moreover, there exists some graph for which the second inequality is strict (see \cite{Schrijver1979}). Given these properties, $\vartheta'(G)$ seems to be a good candidate for upper bounding the Shannon capacity.

More generally, we can find even better approximations for the copositive cone $\mathcal C_n$ by using higher order sum-of-squares polynomials. For each nonnegative integer $r$, define $\mathcal C_n^{(r)}$ to be the set of all symmetric matrices $Q \in \mathcal S_n$ such that
\[
\left(\sum_ix_i^2\right)^rp_Q(x)
\]
is a sum of squares. Then $\mathcal C_n^{(r)}$ is a cone, and
\[
\mathcal P_n \subseteq \mathcal C_n^{(0)} \subseteq \mathcal C_n^{(1)} \subseteq \dots \subseteq \mathcal C_n^{(r)} \subseteq \dots \subseteq \mathcal C_n.
\]
Similar to the Schrijver $\vartheta'$-function, we denote the optimal value of the corresponding problem \eqref{eq:conic} as $\vartheta^{(r)}(G)$.

For higher-order sum-of-squares cones $\mathcal C_n^{(r)}$ where $r>0$, although $\vartheta^{(r)}(G)$ is a tighter upper bound for the independence number $\alpha(G)$ than $\vartheta(G)$ or $\vartheta'(G)$, it is too tight to be an upper bound for the Shannon capacity $\Theta(G)$. For instance, for the pentagon graph $C_5$, if $r>0$,
\[
\alpha(C_5)=\vartheta^{(r)}(C_5)=2<\Theta(C_5)=\vartheta(C_5)=\vartheta'(C_5)=\sqrt 5.
\]
Therefore, to obtain an upper bound for the Shannon capacity from cones $\mathcal C_n^{(r)}$, we have to add extra constraints in the problem \eqref{eq:conic} to restrict these cones. Whatever the exact form of constraints is, we can still analyze the restricted problem as a special case of \eqref{eq:conic} as long as these constraints define a cone.

In the next, we will assume $\mathcal A_n$ in the above problem \eqref{eq:conic} to be an arbitrary cone satisfying $\mathcal A_n \subseteq \mathcal C_n$, and the optimal value will be called $f(G)$. To ensure that $f(G)$ is still an upper bound for the Shannon capacity $\Theta(G)$, we will look at the key property of the semidefinite cone $\mathcal P_n$ used by the Lov\'asz number $\vartheta(G)$ that guarantees $\Theta(G) \leq \vartheta(G)$ and then try to enforce the same property on the cone $\mathcal A_n$ in \eqref{eq:conic}.

\section{Product Property and Upper Bounds for the Shannon Capacity}\label{sec:product}

One fundamental property\footnote{In fact, the equality holds in \eqref{eq:lovaszproduct}, but the reverse direction is not relevant for our purpose.} of the Lov\'asz number is
\begin{equation}\label{eq:lovaszproduct}
\vartheta(G \boxtimes H) \leq \vartheta(G)\vartheta(H)
\end{equation}
for any graphs $G$ and $H$, which immediately implies that
\[
\sqrt[k]{\alpha(G^k)} \leq \sqrt[k]{\vartheta(G^k)} \leq \vartheta(G)
\]
for all positive integers $k$, and thus
\[
\Theta(G)=\sup_k\sqrt[k]{\alpha(G^k)} \leq \vartheta(G).
\]

The above argument can also be applied to the graph function $f(G)$ defined as the optimal value of \eqref{eq:conic}. Since $\alpha(G) \leq f(G)$, as long as $f(G)$ satisfies the similar inequality
\begin{equation}\label{eq:product}
f(G \boxtimes H) \leq f(G)f(H),
\end{equation}
$f(G)$ will also be an upper bound for the Shannon capacity $\Theta(G)$. To find out what leads to the inequality \eqref{eq:product}, we need to generalize the proof for the property \eqref{eq:lovaszproduct} of the Lov\'asz number, which itself is a special case of the general product rules in semidefinite programming \cite{MS2007}.

Consider two graphs $G$ of $n$ vertices and $H$ of $m$ vertices. Assume $(\lambda',Y')$ and $(\lambda'',Y'')$ are the optimal solutions to the problem \eqref{eq:conic} for graph $G$ and $H$, respectively. Let $Y=Y' \otimes Y''$, i.e., the Kronecker product of $Y'$ and $Y''$, which is an $nm \times nm$ matrix given by
\[
Y=\begin{pmatrix}
Y'_{11}Y'' & \cdots & Y'_{1n}Y'' \\
\vdots & \ddots & \vdots \\
Y'_{n1}Y'' & \cdots & Y'_{nn}Y''
\end{pmatrix}.
\]
If we index the rows of $Y$ by pairs $(i,j)$ and the columns by pairs $(k,l)$, the above definition can be rewritten as
\[
Y_{(i,j)(k,l)}=Y'_{ik}Y''_{jl}.
\]
When $(i,j) \not\sim_{G \boxtimes H} (k,l)$, by definition either $i \not\sim_G k$ or $j \not\sim_H l$, which implies either $Y'_{ik}=0$ or $Y''_{jl}=0$ and thus $Y_{(i,j)(k,l)}=0$. Since all the diagonal entries of $Y$ equal to $\lambda'\lambda''$, if we can show that $Y-J_{nm} \in \mathcal A_{nm}$, $(\lambda'\lambda'',Y)$ will be a feasible solution to the problem \eqref{eq:conic} for the product graph $G \boxtimes H$. In this case, we have
\[
f(G \boxtimes H) \leq \lambda'\lambda''=f(G)f(H).
\]

Let
\[
Q=Y'-J_n, \quad R=Y''-J_m.
\]
Then $Q \in \mathcal A_n$, $R \in \mathcal A_m$. The only missing part that remains to be shown is
\[
Y-J_{nm}=Y' \otimes Y''-J_{nm}=(Q+J_n) \otimes (R+J_m)-J_{nm} \in \mathcal A_{nm},
\]
which will be encapsulated into the following definition:

\begin{definition}
Given two symmetric matrices $Q \in \mathcal S_n$, $R \in \mathcal S_m$, define
\[
Q \odot R=(Q+J_n) \otimes (R+J_m)-J_{nm}.
\]
A sequence of cones $\mathcal A_n \subseteq \mathcal S_n$ is said to have the \emph{product property} if for any matrices $Q \in \mathcal A_n$, $R \in \mathcal A_m$, we have $Q \odot R \in \mathcal A_{nm}$.
\end{definition}

Based on this definition, the above argument can be summarized as follows:
\begin{theorem}\label{thm:bound}
If the cones $\mathcal A_n$ in problem \eqref{eq:conic} satisfy $\mathcal A_n \subseteq \mathcal C_n$ and the product property, then $\Theta(G) \leq f(G)$ for any graph $G$.
\end{theorem}

As an example, we check that the product property holds for semidefinite cones $\mathcal P_n$ in the Lov\'asz number. Assume matrices $Q \in \mathcal P_n$, $R \in \mathcal P_m$. Then the matrix
\[
Q \odot R=(Q+J_n) \otimes (R+J_m)-J_{nm}=Q \otimes R+Q \otimes J_m+J_n \otimes R
\]
is also positive semidefinite, because the Kronecker product of two positive semidefinite matrices is still positive semidefinite. Therefore, Theorem~\ref{thm:bound} implies that the Lov\'asz number $\vartheta(G) \geq \Theta(G)$.

The product property is a sufficient condition for the functional inequality \eqref{eq:product} and further for being an upper bound for the Shannon capacity. However, neither the product property nor the inequality \eqref{eq:product} is necessary for being the upper bound. In any case, from the proof of Theorem~\ref{thm:bound}, one can see that the product property is the most natural condition to guarantee $\Theta(G) \leq f(G)$. In the next, we will study the product property holds for what choice of cones $\mathcal A_n$.

\section{Optimality of the Lov\'asz Number}\label{sec:optimal}

In the previous section, we have stated the product property, the condition for our new function $f(G)$ to be an upper bound for the Shannon capacity. At the same time, we do not want $f(G)$ to be much larger than the Lov\'asz number $\vartheta(G)$ for the same graph $G$. Note that the Lov\'asz number satisfies the following sandwich inequality:
\[
\alpha(G) \leq \vartheta(G) \leq \chi(\bar G),
\]
where $\chi(\bar G)$ is the chromatic number for the complement graph of $G$. Choose $G=\bar K_2$, the edgeless graph of two vertices, then
\[
2=\alpha(\bar K_2) \leq \vartheta(\bar K_2) \leq \chi(K_2)=2.
\]
If the new function $f(G)$ satisfies the similar sandwich inequality, we must have $f(\bar K_2)=2$, which means that the matrix
\[
\Lambda=\begin{pmatrix}
1 & -1 \\
-1 & 1
\end{pmatrix} \in \mathcal A_2.
\]

We want to find a sequence of cones $\mathcal A_n$ satisfying all the above desired conditions. However, it turns out that the only possible cones $\mathcal A_n$ must be subsets of the corresponding semidefinite cones $\mathcal P_n$, and consequently the obtained upper bound $f(G)$ would be at least the Lov\'asz number.

\begin{theorem}\label{thm:optimal}
Suppose a sequence of cones $\mathcal A_n$ satisfies the following properties:
\begin{enumerate}
\item $\mathcal A_n \subseteq \mathcal C_n$ for all $n$.
\item The matrix $\Lambda \in \mathcal A_2$.
\item The sequence $\mathcal A_n$ has the product property.
\end{enumerate}
Then we must have $\mathcal A_n \subseteq \mathcal P_n$ for all $n$.
\end{theorem}
\begin{myproof}
We prove by contradiction. Suppose there is a matrix $A \in \mathcal A_n$ that is not positive semidefinite and $v \in \RR^n$ is a vector such that $v^TAv<0$.

The first step is to construct a matrix $B \in \mathcal A_m$ with $m=2n$ and a vector $w \in \RR^m$ satisfying
\[
w^TBw<0, \quad \sum_iw_i=0.
\]
For any $k>0$, let $B=\Gamma \odot (kA)$, then by the cone property and the product property,
\[
B=\begin{pmatrix}
2kA+J_n & -J_n \\
-J_n & 2kA+J_n
\end{pmatrix} \in \mathcal A_m.
\]
If we let
\[
w=\begin{pmatrix}
v \\
-v
\end{pmatrix},
\]
then
\[
w^TBw=4kv^TAv+4v^TJ_nv.
\]
In the above argument, we can choose $k$ with
\[
k>-\frac{v^TJ_nv}{v^TAv},
\]
and now the matrix $B$ and the vector $w$ will have all the desired properties.

Next, we are going to construct another matrix $C \in \mathcal A_{2m}$ which is not copositive. Define
\[
x=\max(w,0), \quad y=\max(-w,0).
\]
Then $x,y \geq 0$ and $w=x-y$. For any $k>1$, by the product property again, the matrix
\[
C=(k\Gamma) \odot B=\begin{pmatrix}
(k+1)B+kJ_m & -(k-1)B-kJ_m \\
-(k-1)B-kJ_m & (k+1)B+kJ_m
\end{pmatrix} \in \mathcal A_{2m}.
\]
Consider
\begin{multline*}
\begin{pmatrix}
x^T & y^T
\end{pmatrix}C\begin{pmatrix}
x \\
y
\end{pmatrix}=(k+1)(x^TBx+y^TBy)-2(k-1)x^TBy \\
+k(x^TJ_mx+y^TJ_my)-2kx^TJ_my,
\end{multline*}
in which the second part
\[
k(x^TJ_mx+y^TJ_my)-2kx^TJ_my=k(x-y)^TJ_m(x-y)=k\left(\sum_iw_i\right)^2=0.
\]
On the other hand, for sufficiently large $k$,
\[
w^TBw=(x-y)^TB(x-y)=x^TBx+y^TBy-2x^TBy<0
\]
implies
\[
x^TBx+y^TBy<2\frac{k-1}{k+1}x^TBy.
\]
In this case,
\[
\begin{pmatrix}
x^T & y^T
\end{pmatrix}C\begin{pmatrix}
x \\
y
\end{pmatrix}<0.
\]
Now we have exhibited a matrix $C \in \mathcal A_{2m}$ and $C$ is not copositive, which is a contradiction.
\end{myproof}

Theorem~\ref{thm:optimal} tells us that either the cones do not have the product property or the resulting function $f(G) \geq \vartheta(G)$. In other words, it is impossible to derive an upper bound for the Shannon capacity that is better than the Lov\'asz number by enforcing the product property on cones $\mathcal A_n$.

For the Schrijver $\vartheta'$-function, the corresponding cones $\mathcal C_n^{(0)}$ satisfy the first and second condition of Theorem~\ref{thm:optimal} but not the conclusion $\mathcal C_n^{(0)} \subseteq \mathcal P_n$. Therefore, by Theorem~\ref{thm:optimal}, the cones $\mathcal C_n^{(0)}$ do not have the product property. Although not having the product property for $\mathcal C_n^{(0)}$ does not directly imply that the Schrijver $\vartheta'$-function is not an upper bound for the Shannon capacity, it strongly suggests such negative result. In fact, it is quite difficult to disprove that the Schrijver $\vartheta'$-function is an upper bound, because at least for graphs $G$ of moderate size the two values $\vartheta'(G)$ and $\vartheta(G)$ are very close to each other. In order to prove that $\vartheta'(G)$ is not an upper bound, we have to find some sufficiently large $k$ and show that
\[
\vartheta'(G)<\sqrt[k]{\alpha(G^k)} \leq \vartheta(G),
\]
which is extremely hard even if $G$ contains only a few vertices. We believe that the Schrijver $\vartheta'$-function is not an upper bound for the Shannon capacity due to its lack of the product property, but whether this is actually true or not remains open.

\bibliography{SOS}
\bibliographystyle{spmpsci}
\end{document}